\documentclass[12pt, english]{article}
\usepackage[all]{xy}
\usepackage{varioref,babel}
\usepackage{enumerate,a4,amsfonts,theorem,amsmath,amsfonts,amssymb,amscd}
\usepackage{amsxtra,eucal}
\usepackage[latin1]{inputenc}
\usepackage{tabularx}

\begin{document}
\pagestyle{plain}
\parindent0mm

\title{Hermitian categories, extension of scalars and systems of sesquilinear forms}
\author{Eva Bayer-Fluckiger and Daniel Arnold Moldovan\footnote{Partially supported by the
Swiss National Science Fundation, grant 200020-109174/1} \\ \'Ecole Polytechnique
F\'ed\'erale de
Lausanne, Switzerland}
\noindent
\maketitle

\textbf{Abstract }  In this paper we define a notion of Witt group for sesquilinear forms in hermitian categories, which in turn provides a notion of Witt group for sesquilinear forms over rings
with involution. We also study the extension of scalars for $K$-linear hermitian categories, where $K$ is a field of characteristic $\neq 2$.
We finally extend several results concerning sesquilinear forms to the setting of systems of 
such forms.

\bigskip

\textbf{Mathematics Subject Classification (2000) } 11E39, 11E81. \medskip \\

\textbf{Keywords } sesquilinear forms, hermitian forms, systems of sesquilinear forms, hermitian categories, $K$-linear categories, scalar extension, Witt group.

\medskip
\bigskip

{\bf \par Introduction}
\bigskip \\
There exists a classical notion of Witt group for hermitian forms over rings with involution (see for instance [6] and [8]). 
However, there is no analogous notion for sesquilinear forms. In this paper we define, based on an equivalence of categories proven in [5], 
a notion of Witt group for sesquilinear forms in hermitian categories. In particular we obtain a notion of Witt group for sesquilinear forms over a ring with involution.

\medskip

We then study the extension of scalars for $K$-linear hermitian categories, where $K$ is a field of characteristic $\neq 2$, and we prove its injectivity in the case of an extension
of odd degree and a hermitian category in which all idempotents split.

\medskip

We also introduce the notion of system of sesquilinear forms over a ring with involution and we generalize several results proven in [5] to systems of sesquilinear forms.

\bigskip
\bigskip

\newpage

{\bf \S 1. Sesquilinear and hermitian forms over rings with involution}
\medskip \\

Let $A$ be a ring. An {\it involution} on $A$ is by definition an additive map  $\sigma : A \to A$ 
such that $\sigma(ab) = \sigma(b) \sigma(a)$ for all $a,b \in A$ and $\sigma^2$ is the identity. 
Let $V$ be a right $A$--module of finite type. 
A {\it sesquilinear form} over $(A,\sigma)$ is a biadditive map
$s : V \times V \to A$ satisfying the condition $s(ax,by) = \sigma(a) s(x,y) b$ for all $x,y \in V$ and all $a, b \in A$. The {\it orthogonal sum} of two sesquilinear forms $(V,s)$ and $(V',s')$
is by definition the form $(V \oplus V', s \oplus s')$ defined by 
$$(s \oplus s')(x \oplus x', y \oplus y') =
s(x,y) + s'(x',y')$$
for all $x, y \in V$ and $x', y' \in V'$. Two sesquilinear forms $(V,s)$ and $(V', s')$ are called \textit{isometric} if
there exists an isomorphism of $A$-modules $f: V \tilde{\rightarrow} V'$ such that $s'(f(x), f(y))=s(x,y)$ for all $x, y \in V$.
\medskip \\
Let $V^* = {\rm Hom}_A(V,A)$. Then $V^*$ has a structure of right $A$--module given by
$f(x) \cdot a = \sigma(a) f(x)$ for all $a \in A$, $f \in V^*$ and $a \in A$. We say that $V$ is
{\it reflexive} if the homomorphism of right $A$-modules $e_V : V \to V^{**}$ defined by $e_V(x)(f) = \sigma(f(x))$ for all $x \in V$ and $f \in V^*$
is bijective. 
\medskip \\
A sesquilinear form $(V,s)$ over $(A, \sigma)$ induces two homomorphisms of left $A$-modules $V \to V^*$, called its {\it left}, respectively  
{\it right adjoint}, namely $s_l : V \to V^*$ defined by $s_l(x)(y) = s(x,y)$ and
$s_r : V \to V^*$ given by $s_r(x)(y) = \sigma(s(y,x))$ for all $x, y \in V$. We observe that $s_r=s_l^* e_V$.
\medskip \\
Let $\epsilon=\pm 1$. A sesquilinear form $(V,s)$ over $(A, \sigma)$ is called \textit{$\epsilon$-hermitian} if $V$ is a projective $A$-module and 
$\sigma(s(x,y))=\epsilon s(y,x)$ for all $x, y \in V$, i.e. $s_l=\epsilon s_r$. A $1$-hermitian form is also called a hermitian form. 
An $\epsilon$-hermitian form $(V,s)$ is called \textit{unimodular} if $s_l$ (or equivalently $s_r$) is bijective. 
There exists a classical notion of Witt group for unimodular $\epsilon$-hermitian forms over $(A, \sigma)$ (see e.g. [6]). Denote by ${\rm Gr}^{\epsilon}(A, \sigma)$ the Grothendieck group
of isometry classes of unimodular $\epsilon$-hermitian forms over $(A, \sigma)$, with respect to the orthogonal sum. A unimodular $\epsilon$-hermitian form
over $(A, \sigma)$ is called \textit{hyperbolic} if it is isometric to a form $\mathbb{H}(V)$, where $V$ is a finitely generated projective right $A$-module and 
$$\mathbb{H}(V): V \otimes V^* \rightarrow V^*\oplus V^{**}$$
$$x+y \mapsto y+ \epsilon e_V(x), ~\forall x \in V, y \in V^*.$$
The quotient of ${\rm Gr}^{\epsilon}(A, \sigma)$ by the subgroup generated by the unimodular $\epsilon$-hermitian forms is called the \textit{Witt group of unimodular $\epsilon$-hermitian forms} 
over $(A, \sigma)$ and is denoted by $W^{\epsilon}(A, \sigma)$.
\medskip \\
Let us denote by $\mathcal{S}(A, \sigma)$ (~$\mathcal{H}^{\epsilon}(A, \sigma)$~) the category of sesquilinear (respectively unimodular $\epsilon$-hermitian) forms over $(A, \sigma)$. 
The morphisms of these categories are isometries. For simplicity let $\mathcal{H}(A, \sigma)=\mathcal{H}^1(A, \sigma)$.

\bigskip
\bigskip

{\bf \S 2. Hermitian categories} 
\medskip \\

The aim of this section is to recall some basic notions about hermitian categories as presented in [8] (see also [6], [7]). 

\bigskip

{\bf \S 2.1. Preliminaries}

\bigskip

Let $\cal C$ be an additive category. Let
${}^* : {\cal C} \to {\cal C}$ be a duality functor, i.e. an additive contravariant functor with
a natural isomorphism $(E_C) _{C \in {\cal C}} : {\rm id} \to {}^{**}$ such that $E^*_C E_{C^*} 
= {\rm id}_{C^*}$ for all $C \in {\cal C}$. An additive category with a duality functor is called a \textit{hermitian category}. 
A {\it sesquilinear form} in the category ${\cal C}$ is a pair $(C, s)$, where $C$ is an object of ${\cal C}$
and $s: C \rightarrow C^*$. A sesquilinear form $(C, s)$ is called {\it unimodular} if $s$ is an isomorphism. Let $\epsilon=\pm 1$. An {\it $\epsilon$-hermitian form}
in the category ${\cal C}$ is a sesquilinear form $(C, s)$ such that $s= \epsilon s^* E_C$. A $1$-hermitian form is also called a hermitian form. Orthogonal sums of 
forms are defined in the obvious way. Let $(C,s)$ and $(C',s')$ be two sesquilinear forms in ${\cal C}$. We say that these forms
are {\it isometric} if there exists an isomorphism $f : C \tilde{\rightarrow} C'$ in the category $\cal C$
such that $s = f^* s' f$. 
\medskip \\
Denote by $\mathcal{H}^{\epsilon}(\cal C)$ the category of unimodular $\epsilon$-hermitian forms in the category $\cal C$. 
The morphisms are isometries. For simplicity let $\mathcal{H}(\mathcal{C})=\mathcal{H}^1(\mathcal{C})$.
\medskip \\
The hyperbolic unimodular $\epsilon$-hermitian forms in $\mathcal{C}$ are the forms isometric to $\mathbb{H}_{Q}$, $Q \in \mathcal{C}$, given by
$$\mathbb{H}_{Q}=\epsilon E_Q \oplus {\rm id}_{Q^*}: Q \oplus Q^* \rightarrow Q^{**} \oplus Q^* \simeq (Q \oplus Q^{*})^*.$$
The quotient of the Grothendieck group of isometry classes of unimodular $\epsilon$-hermitian forms in $\mathcal{C}$ (with respect to the orthogonal sum) 
by the subgroup generated by the hyperbolic forms is called the Witt group of 
unimodular $\epsilon$-hermitian forms in $\mathcal{C}$ and is denoted by $W^{\epsilon}(\mathcal{C})$. For simplicity set $W(\mathcal{C})=W^1(\mathcal{C})$.
\medskip \\
We observe that if we take $\mathcal{C}$ to be the category of reflexive right $A$-modules of finite type, then the notion of sesquilinear form coincides with the one defined in the 
preceding section. Analogously, if we take $\mathcal{C}$ to be the category of projective right $A$-modules of finite type, then the notion of hermitian form coincides with the one
defined in the preceding section. 
\medskip \\
Let $(\mathcal{M}, *)$ and $(\mathcal{M}', *)$ be two hermitian categories. A \textit{duality preserving functor} from $(\mathcal{M}, *)$ to $(\mathcal{M}', *)$ is an additive functor
$F: \mathcal{M} \rightarrow \mathcal{M}'$ together with a natural isomorphism $i=(i_M)_{M \in \mathcal{M}}:F^* \rightarrow ~^*F$. This means that for any $M \in \mathcal{M}$ there exists a
natural isomorphism $i_M:F(M^*) \tilde{\rightarrow} F(M)^*$ such that for all $N \in \mathcal{M}$ and $f \in {\rm Hom}_{\mathcal{M}}(M,N)$ the
following diagram commutes:

$$\xymatrix@!C{
F(N^*) \ar[r]^{F(f^*)} \ar[d]^{i_N} & F(M^*) \ar[d]^{i_M}\\
F(N)^* \ar[r]^{F(f)^*} & F(M)^*
}$$

\medskip

For $\lambda=\pm 1$, a duality preserving functor $F$ is called
\textit{$\lambda$-hermitian} if $i_{M^*} F(e_M)=\lambda i^*_M e_{F(M)}$ for all $M \in \mathcal{M}$. Let $\epsilon=\pm 1$. We recall from [6], pp. 80-81 that
a $\lambda$-hermitian functor $F: \mathcal{M} \rightarrow \mathcal{M}'$ induces a functor
$${\rm Herm}^{\epsilon}(F): \mathcal{H}^{\epsilon}(\mathcal{M}) \rightarrow \mathcal{H}^{\epsilon \lambda}(\mathcal{M}')$$
$$(M, h) \mapsto (F(M), i_M F(h)),$$
which preserves orthogonal sums and hyperbolicity. Therefore it induces a homomorphism between the corresponding Witt groups:
$${\rm \overline{Herm}}^{\epsilon}(F): W^{\epsilon}(\mathcal{M}) \rightarrow W^{\epsilon \lambda}(\mathcal{M}').$$

\smallskip

If $F$ is fully faithful, then ${\rm Herm}^{\epsilon}(F)$ is also fully faithful. Moreover, if $F$ is an equivalence of categories, then ${\rm Herm}^{\epsilon}(F)$ 
is also an equivalence of categories and the induced group homomorphism ${\rm \overline{Herm}}^{\epsilon}(F)$ is bijective. 

\bigskip
\bigskip

{\bf \S 2.2. Transfer into the endomorphism ring}
\medskip \\
The aim of this subsection is to introduce the method of \textit{transfer into the endomorphism ring}, which allows us to pass from the abstract setting of hermitian categories to that of a ring with involution,
which is more concrete. This method will be extensively applied in section 4.
\medskip \\
Let $\mathcal{M}$ be a hermitian category and $M$ be an object of $\mathcal{M}$, on which we suppose that there exists a unimodular $\epsilon_0$-hermitian form $h_0$ for a certain $\epsilon_0=\pm 1$. 
Denote by $E$ the endomorphism ring of $M$. According to [7], lemma 1.2, the form $(M,h)$ induces on $E$ an involution $\sigma$, defined by
$\sigma(f)=h_0^{-1} f^* h_0$ for all $f \in E$. 
\medskip \\
We say that an \textit{idempotent $e \in E$ splits} if there exist an object $M' \in \mathcal{M}$ and morphisms 
$i:M' \rightarrow M$, $j: M \rightarrow M'$ such that $ji={\rm id}_{M'}$ and $ij=e$. 
\medskip \\
Denote by $\mathcal{M}\vert_M$ the full subcategory of $\mathcal{M}$ which has as objects all the objects of $\mathcal{M}$ isomorphic to a direct summand
of a finite direct sum of copies of $M$ and by $\mathcal{P}(E)$ the category of finitely generated projective right $E$-modules. We consider the following functor:
$$\mathcal{F}={\rm Hom}(M,-): \mathcal{M}\vert_M \rightarrow \mathcal{P}(E)$$
$$N \mapsto {\rm Hom}(M, N), ~\forall N \in \mathcal{M}\vert_M$$
$$f \mapsto \mathcal{F}(f), ~\forall f \in {\rm Hom}(N,N'), ~\forall N,N' \in \mathcal{M}\vert_M,$$
where for all $g \in {\rm Hom}(M,N)$, $\mathcal{F}(f)(g)=fg$.
In [7], proposition 2.4 it has been proven that this functor is fully faithful and duality preserving with respect to the natural isomorphism 
$i=(i_N)_{N \in \mathcal{M}\vert_M}:\mathcal{F}^* \rightarrow ~^*\mathcal{F}$ given by $i_N(f)=\mathcal{F}(h_0^{-1}f^* e_N)$ for every $N \in \mathcal{M}\vert_M$ and $f \in {\rm Hom}(M,N^*)$. In addition,
if all the idempotents of $\mathcal{M}\vert_N$ split, then $\mathcal{F}$ is an equivalence of categories. By computation we easily see that $\mathcal{F}$ is $\epsilon_0$-hermitian. 

\bigskip
\bigskip

{\bf \S 2.3. Ring extension and genus in hermitian categories}

\medskip
In this subection we introduce the notions of extension of rings and genus in hermitian categories. 
\medskip \\
Let $\mathcal{C}$ be an additive category and $R$ be a commutative ring with unity. By extension of scalars from
$\mathbb{Z}$ to $R$ we obtain a new category, denoted by $\mathcal{C} \otimes_{\mathbb{Z}} R$ and called the \textit{extension of the category $\mathcal{C}$ to the ring $R$}. 
Its objects are the same as those of $\mathcal{C}$ and for two such objects $M$ and $N$ set 
$${\rm Hom}_{\mathcal{C} \otimes_{\mathbb{Z}} R} (M, N)={\rm Hom}_{\mathcal{C}}(M, N) \otimes_{\mathbb{Z}} R.$$
It is straightforward to check that in this way we obtain an additive category. If, in addition, there is a duality functor $*$ in the category $\mathcal{C}$, then
$\mathcal{C} \otimes_{\mathbb{Z}} R$ becomes a hermitian category by setting $(f \otimes a)^*=f ^* \otimes a$ for all $M, N \in \mathcal{C}$, $f \in {\rm Hom}_{\mathcal{C}}(M, N)$ and 
$a \in R$. \medskip \\
For each prime spot $p$ of $\mathbb{Q}$ we denote by $\mathcal{C}_p$ the category $\mathcal{C} \otimes_{\mathbb{Z}} \mathbb{Z}_p$. \medskip \\
Let $\mathcal{C}$ be an additive category. Two sesquilinear forms $s: M \rightarrow M^*$ and $t:N \rightarrow N^*$ in the category $\mathcal{C}$ are said to be in the same \textit{genus} if for all primes $p$ of $\mathbb{Q}$ there exist
$F_p \in {\rm Hom}_{\mathcal{C}}(M, N)_p$ and $G_p \in {\rm Hom}_{\mathcal{C}}(N, M)_p$ such that $F_p G_p$ and $G_p F_p$ are equal to the identity and in addition $F_p^* (t \otimes 1) F_p=s \otimes 1$. This is equivalent to
saying that the extended sesquilinear forms $(M, s \otimes 1)$ and $(N, t \otimes 1)$ become isometric in the extended category $\mathcal{C}_p$ for every $p$.

\bigskip
\bigskip

{\bf \S 2.4. $K$-linear hermitian categories and scalar extension}
\medskip \\
The aim of this section is to give an introduction to the theory of $K$-linear hermitian categories.
Their consideration is motivated by the idea of defining a notion
of scalar extension in hermitian categories. 
\medskip \\
\textbf{2.4.1. Definition } Let $K$ be a field. A \textit{$K$-linear category} is an additive
category
$\mathcal{M}$ such that for every $M, N \in \mathcal{M}$, the set ${\rm Hom}_{\mathcal{M}}(M, N)$
has a structure
of finite-dimensional $K$-vector space such that the composition of morphisms is $K$-bilinear. An
additive functor $F: \mathcal{M} \rightarrow \mathcal{N}$ between two $K$-linear
categories is called \textit{$K$-linear} if for any $M, N \in \mathcal{M}$, the induced map
$$F_{M, N}: {\rm Hom}_{\mathcal{M}}(M, N) \rightarrow {\rm
Hom}_{\mathcal{N}}(F(M), F(N))$$ is
$K$-linear.

\medskip

A \textit{$K$-linear hermitian category} is a $K$-linear category together with a $K$-linear duality
functor.

\medskip

Let $(\mathcal{M}, *)$ be a $K$-linear hermitian category and consider a finite field extension $L$
of $K$. We define the \textit{extension of the category $\mathcal{M}$ to $L$} as being the category
$\mathcal{M}_L$ with the same objects as $\mathcal{M}$ and with the morphisms given by
$${\rm Hom}_{\mathcal{M}_L}(M, N)={\rm Hom}_{\mathcal{M}}(M, N) \otimes_K L$$
for all $M, N \in \mathcal{M}$. It is clear that the category $\mathcal{M}_L$ is $L$-linear.

\medskip

If in $\mathcal{M}$ there is a $K$-linear duality functor $*$ then we can define a duality functor in $\mathcal{M}_L$ in the following way:
on objects it coincides with the duality functor of $\mathcal{M}$ and for morphisms set $(f \otimes a)^*=f^*
\otimes a$ for all $f \in {\rm Hom}_{\mathcal{M}}(M, N)$, $M, N \in \mathcal{M}$ and $a \in L$ and then
extend by additivity.

\medskip

The \textit{scalar extension functor} from $\mathcal{M}$ to $\mathcal{M}_L$ is defined by
$$\mathcal{R}_{L/K}: \mathcal{M} \rightarrow \mathcal{M}_L$$
$$M \rightarrow M, ~\forall M \in \mathcal{M}$$
$$f \mapsto f \otimes 1, ~\forall f \in {\rm Hom}(M,N).$$




\medskip 

It is straightforward to check that the functor $\mathcal{R}_{L/K}$ is $1$-hermitian. As we have seen in \S 2.1, it induces for every $\epsilon=\pm 1$ a group homomorphism
$${\rm \overline{Herm}}^{\epsilon}(\mathcal{R}_{L/K}):W^{\epsilon}(\mathcal{M}) \rightarrow W^{\epsilon}(\mathcal{M}_L),$$
called \textit{restriction map}. 

\bigskip 
\bigskip

{\bf \S 3. Sesquilinear forms over rings with involution and hermitian categories}

\medskip 

In this section we prove that the category of sesquilinear forms over a ring with involution is equivalent with the category of unimodular hermitian forms in a suitably constructed 
hermitian category. We also define a notion of Witt group for sesquilinear forms in hermitian categories. 
In particular we obtain a notion of Witt group for sesquilinear forms over
rings with involution, which generalizes the analogous notion for unimodular hermitian forms. 

\medskip

{\bf \S 3.1. An equivalence of categories}

\medskip

Let $\mathcal{M}$ be an additive category. On the model of [5], \S 3, we construct the \textit{category of double arrows of $\mathcal{M}$}, denoted by $\mathcal{M}^{(2)}$.
Its objects are of the form $(M,N,f,g)$, where
$M, N \in \mathcal{M}$ and $f,g \in {\rm Hom}(M,N)$. A morphism from $(M,N,f,g)$ to $(M',N',f',g')$ is a pair $(\phi, \psi)$, where
$\phi \in {\rm Hom}(M, M')$ and $\psi \in {\rm Hom}(N,N')$ satisfying the conditions $\psi f=f' \phi$ and $\psi g=g' \phi$. Then $\mathcal{M}^{(2)}$ is obviously an additive category. If in addition there is a
duality functor $*$ in $\mathcal{M}$, then $\mathcal{M}^{(2)}$ becomes a hermitian category by defining the dual of $(M,N,f,g)$ as being $(N^*, M^*, g^*, f^*)$ and setting
$E_{(M,N,f,g)}=(e_M, e_N)$. 

\medskip 

We define a functor $F: \mathcal{S}(\mathcal{M}) \rightarrow \mathcal{H}(\mathcal{M}^{(2)})$ as follows. Let $(M,s)$ be a sesquilinear form in $\mathcal{M}$. Then
$(M, M^*, s, s^* e_M)$ is an object of $\mathcal{M}^{(2)}$ and it is easy to check that $(e_M, {\rm id}_{M^*})$ defines a unimodular hermitian form on it. Set
$F(M,s)=((M, M^*, s, s^* e_M), (e_M, {\rm id}_{M^*}))$. For an isometry $\phi:(M,s) \rightarrow (M',s')$ of sesquilinear forms in $\mathcal{M}$ set
$F(\phi)=(\phi, \phi^{*-1})$. 

\medskip

We also define a functor $G$ in the opposite direction: for every object $(M,N,f,g)$ of $\mathcal{M}^{(2)}$ and any unimodular hermitian form $(\xi_1, \xi_2)$ on it set
$G((M,N,f,g), (\xi_1,\xi_2))=(M, \xi_2 f)$. For a morphism $(\lambda_1, \lambda_2)$ between two unimodular hermitian forms in $\mathcal{M}^{(2)}$ define $G(\lambda_1, \lambda_2)=\lambda_1$.

\medskip

\textbf{3.1.1. Theorem } \textit{The functors $F$ and $G$ realize an equivalence between the categories $\mathcal{S}(\mathcal{M})$ and $\mathcal{H}(\mathcal{M}^{(2)})$.}

\medskip

\textbf{Proof. } The proof is completely analogous to the one of [5], theorem 4.1. 

\bigskip
\bigskip

{\bf \S 3.2. Hyberbolic sesquilinear forms}

\medskip

Using the functor $G$ we can define a notion of hyperbolicity for sesquilinear forms in $\mathcal{M}$, which will in turn lead to a notion of Witt group. 

\medskip 

\textbf{3.2.1. Definition } A sesquilinear form in $\mathcal{M}$ is called \textit{hyperbolic} if it is isometric to a form $G(P,s)$, where
$(P, s) \in \mathcal{H}(\mathcal{M}^{(2)})$ is hyperbolic.

\medskip

Let $(A, \sigma)$ be a ring with involution. Taking $\mathcal{M}$ to be the category of reflexive right $A$-modules of finite type we obtain a notion of hyperbolicity for sesquilinear
forms over $(A, \sigma)$.

\medskip 

We have the following explicit characterization of hyperbolic sesquilinear forms over $(A, \sigma)$:
\medskip \\
\textbf{3.2.2 Proposition } \textit{A sesquilinear form $(V, s)$ over $(A, \sigma)$ is hyperbolic if and only if there exist two reflexive right $A$-modules of finite type
$M$ and $N$ and two $A$-linear homomorphisms $f, g: M \rightarrow N$ such that $(V,s)$ is isometric to the form $(M \oplus N^*, \widetilde{H}_{(M,N,f,g)})$,
given by
$$\widetilde{H}_{(M, N, f, g)}: (M \oplus N^*) \times (M \oplus N^*) \rightarrow A$$
$$(x_1+y_1, x_2+y_2) \mapsto y_1(g(x_2))+\sigma(y_2(f(x_1))), ~\forall x_1, x_2 \in M, ~\forall y_1, y_2 \in N^*.$$}
\smallskip \\
\textbf{Proof. } For all $Q=(M, N, f, g) \in \mathcal{M}^{(2)}$ the hyperbolic unimodular hermitian form $\mathbb{H}_Q$ in $\mathcal{M}^{(2)}$ is defined by 
$$\mathbb{H}_{Q}=E_Q \oplus {\rm id}_{Q^*}: Q \oplus Q^* \rightarrow Q^{**} \oplus Q^{*} \simeq (Q \oplus Q^*)^*,$$
hence $G(\mathbb{H}_Q)$ is given by
$$G(\mathbb{H}_Q): M \oplus N^* \rightarrow M^* \oplus N^{**}$$
$$x+y \mapsto g^*(y)+e_N(f(x)), ~\forall x \in M, ~\forall y \in N^*.$$
It is easy to identify this form with $\widetilde{H}_Q$. $\square$
\medskip \\
Particularizing the last proposition, we obtain:
\bigskip \\
\textbf{3.2.3. Proposition } \textit{
\begin{enumerate}[a)]
 \item The hyperbolic hermitian forms over $(A, \sigma)$ are given, up to isometry, by the forms $(M \oplus N^*, \widetilde{H}_{(M, N, f, f)})$, where $M$ and $N$ are two
reflexive right $A$-modules of finite type and $f: M \rightarrow N$ is an $A$-linear homomorphism.
 \item The hyperbolic unimodular sesquilinear forms over $(A, \sigma)$ are given, up to isometry, by $(M \oplus N^*, \widetilde{H}_{(M, N, f, g)})$, where $M$
and $N$ are two reflexive
right $A$-modules of finite type and $f, g: M \rightarrow N$ are two $A$-linear isomorphisms.
 \item The hyperbolic unimodular hermitian forms over $(A, \sigma)$ are given, up to isometry, by $(M \oplus M^*, \mathbb{H}_M)$, where $M$ is a reflexive right
$A$-module of finite type and
$$\mathbb{H}_M: (M \oplus M^*) \times (M \oplus M^*) \rightarrow A$$
$$(x_1+y_1, x_2+y_2) \mapsto y_1(x_2)+\sigma(y_2(x_1)), ~\forall x_1, y_1 \in M, ~\forall y_1, y_2 \in M^*.$$
\end{enumerate}
}
\medskip 
\textbf{Proof. } It is straightforward to check that the form $(M \oplus N^*, \widetilde{H}_{(M, N, f, g)})$ is hermitian if and only if $f=g$ and that it is unimodular
if and only if $f$ and $g$ are isomorphisms. The point c) follows from a) and b) and from the fact that the forms $(M \oplus N^*, \widetilde{H}_{(M, N, f, f)})$
and $(M \oplus M^*, \mathbb{H}_M)$ are isometric via the isomorphism of $A$-modules ${\rm id}_{M} \oplus f^*: M \oplus N^* \tilde{\rightarrow} M \oplus M^*$. $\square$
\bigskip \\
\textbf{3.2.4. Definition } Let ${\rm Gr}_S(\mathcal{M})$ be the Grothendieck group of isometry classes of sesquilinear forms in $\mathcal{M}$, with respect to the orthogonal sum. 
It is easy to check that the isometry classes
of hyperbolic sesquilinear forms form a subgroup of ${\rm Gr}_S(\mathcal{M})$, which we denote by $\mathbb{H}_S(\mathcal{M})$. The \textit{Witt group of
sesquilinear forms} in the category $\mathcal{M}$ is defined to be the quotient
$$W_S(\mathcal{M})={\rm Gr}_S(\mathcal{M})/\mathbb{H}_S(\mathcal{M}).$$

\smallskip

Taking $\mathcal{M}$ to be the category of reflexive right $A$-modules of finite type we obtain a notion of Witt group for sesquilinear forms over $(A, \sigma)$. 
We observe that in the case of unimodular hermitian forms our definition coincides with the well-known one.

\bigskip
\bigskip

{\bf \S 3.3. Finiteness results concerning the genus of sesquilinear forms}

\medskip

In this subsection we prove a finiteness result concerning the genus of sesquilinear forms, based on the methods of [5], \S 9. For a ring $A$ we denote by $T(A)$ the $\mathbb{Z}$-torsion subgroup of $A$. If $R$ is a ring containing $\mathbb{Z}$, then we say
that $A$ is $R$-\textit{finite} if $A_R=A \otimes_{\mathbb{Z}} R$ is a finitely generated $R$-module and $T(A)$ is finite.
\medskip \\
\textbf{3.3.1. Theorem } \textit{Let $A$ be a ring and $\sigma$ be an involution on $A$. Let $(V,s)$ be a sesquilinear form over $(A, \sigma)$ and assume that ${\rm End}_A(V)$ is $\mathbb{Q}$-finite. Then the genus of $(V, s)$ contains only a finite number of isometry classes of
sesquilinear forms.} \medskip \\
\textbf{Proof. } The functor $F$ defined in \S 3.1 induces a bijection between the genus of $(V,s)$ and the genus of $F(V,s)$. Denote by $E$ the endomorphism ring of $(V, V^*, s_l, s_r)$ in $\mathcal{M}$.
Since ${\rm End}_A(V)$ is $\mathbb{Q}$-finite and $E$ is a subring of ${\rm End}_A(V) \times {\rm End}_A(V^*) \simeq {\rm End}_A(V) \times {\rm End}_A(V)$, $E$ is $\mathbb{Q}$-finite too. 
From [3], theorem 3.4 it follows that the genus of $F(V,s)$ contains only a finite number of isometry classes
of unimodular hermitian forms, which implies that the genus of $(V, s)$ contains only a finite number of isometry classes of sesquilinear forms. 

\bigskip
\bigskip

{\bf \S 4. The restriction map for $K$-linear hermitian categories in odd degree extensions}
\medskip \\
It is well-known that if $K$ is a field of characteristic $\neq 2$ and $L/K$ is an extension of odd degree, then the restriction map $r_{L/K}: W(K) \rightarrow W(L)$ is injective. 
The aim of this section is to prove an analogous result for $K$-linear hermitian categories.
\medskip \\
Let $K$ be a field of characteristic $\neq 2$ and $\mathcal{M}$ be a $K$-linear hermitian category. Let $L$ be a finite extension of $K$, $\mathcal{M}_L$ be the extension of the 
category $\mathcal{M}$ to $L$ and $\mathcal{R}_{L/K}: \mathcal{M} \rightarrow \mathcal{M}_L$ be the scalar extension functor, as defined in \S 2.5. 
\medskip \\
\textbf{4.1. Theorem } \textit{Suppose that all the idempotents of $\mathcal{M}$ split and that the extension $L/K$ is of odd degree. Then for all $\epsilon=\pm 1$ the map
$${\rm \overline{Herm}}^{\epsilon}(\mathcal{R}_{L/K}): W^{\epsilon}(\mathcal{M}) \rightarrow
W^{\epsilon}(\mathcal{M}_L)$$
is injective.}
\medskip \\
This result will follow as an immediate corollary from the following one. 
\newpage

\textbf{4.2. Theorem } \textit{Let $M$ be an object of $\mathcal{M}$ that admits a unimodular hermitian or skew-hermitian form. 
Suppose that all the idempotents of $\mathcal{M}\vert_M$ split and that the extension $L/K$ is of odd degree. Then for all $\epsilon=\pm 1$, the map
$${\rm \overline{Herm}}^{\epsilon}(\mathcal{R}^M_{L/K}): W^{\epsilon}(\mathcal{M}\vert_M) \rightarrow W^{\epsilon}(\mathcal{M}_L\vert_M)$$ is injective.
}

\medskip 

\textbf{Proof of theorem 4.2 } Denote by $h_0$ a unimodular $\epsilon_0$-hermitian form on $M$. 
Since the category $\mathcal{M}$ is $K$-linear, the endomorphism ring $E$ of $M$ in $\mathcal{M}$ has a structure of finite-dimensional $K$-algebra. 
According to [7], lemma 1.2, $h_0$ induces on $E$ a $K$-linear involution, denoted by $\sigma$. 
Clearly the endomorphism ring $E'$ of $M$ in $\mathcal{M}_L$ equals $E \otimes_K L$ and is a finite-dimensional $L$-algebra. We also observe that $(M, h_{0L})$ is a unimodular $\epsilon_0$-hermitian
form in $\mathcal{M}_L$ and that it induces by [7], lemma 1.2 the involution $\sigma \otimes {\rm id}_L$ on $E'$. We have an obvious $1$-hermitian functor of extension of scalars
$$\mathcal{R}_{L/K}^M: \mathcal{M}\vert_M \rightarrow \mathcal{M}_L\vert_M.$$
We recall from \S 2.2 that there are two fully faithful $\epsilon_0$-hermitian functors
$$\mathcal{F}: \mathcal{M}\vert_M \rightarrow \mathcal{P}(E),$$
$$\mathcal{F}': \mathcal{M}_L\vert_M \rightarrow \mathcal{P}(E \otimes_K L).$$
We denote by $R_{L/K}$ the obvious functor of scalar extension 
$$R_{L/K}: \mathcal{P}(E) \rightarrow \mathcal{P}(E \otimes_K L).$$

It is straightforward to prove that the following diagram commutes:
$$\xymatrix@!C{
\mathcal{M}_L \vert_M \ar[r]^-{\mathcal{F}'} & \mathcal{P}(E \otimes_K L)\\
\mathcal{M}\vert_M \ar[r]^-{\mathcal{F}} \ar[u]^{\mathcal{R}^M_{L/K}} & \mathcal{P}(E)
\ar[u]_{R_{L/K}}
}$$

For $\epsilon=\pm 1$ the above diagram induces a commutative diagram on the level of unimodular hermitian forms:

$$\xymatrix@!C{
\mathcal{H}^{\epsilon}(\mathcal{M}_L \vert_M) \ar[r]^-{{\rm Herm}^{\epsilon}(\mathcal{F}')} &
\mathcal{H}^{\epsilon \epsilon_0}(E \otimes_K L, \sigma \otimes {\rm id_L})\\
\mathcal{H}^{\epsilon}(\mathcal{M}\vert_M) \ar[r]^-{{\rm Herm}^{\epsilon}(\mathcal{F})}
\ar[u]^{{\rm Herm}^{\epsilon}(\mathcal{R}^M_{L/K})} & \mathcal{H}^{\epsilon
\epsilon_0}(E, \sigma)
\ar[u]_{{\rm Herm}^{\epsilon \epsilon_0}(R_{L/K})}
}$$

and one on the level of Witt groups:

$$\xymatrix@!C{
W^{\epsilon}(\mathcal{M}_L \vert_M) \ar[r]^-{{\rm \overline{Herm}}^{\epsilon}(\mathcal{F}')} & W^{\epsilon
\epsilon_0}(E \otimes_K L, \sigma \otimes {\rm id}_L)\\
W^{\epsilon}(\mathcal{M}\vert_M) \ar[r]^-{{\rm \overline{Herm}}^{\epsilon}(\mathcal{F})} \ar[u]^{{\rm
\overline{Herm}^{\epsilon}}(\mathcal{R}^M_{L/K})} & W^{\epsilon \epsilon_0}(E, \sigma)
\ar[u]_{{\rm \overline{Herm}^{\epsilon \epsilon_0}}(R_{L/K})}
}$$

\medskip 

Let $(N, g)$ and $(N', g') \in \mathcal{H}^{\epsilon}(\mathcal{M}\vert_M)$ be such that
$${\rm \overline{Herm}}^{\epsilon}(\mathcal{R}^M_{L/K})(N, g)={\rm \overline{Herm}}^{\epsilon}(\mathcal{R}^M_{L/K})(N', g')$$ in
$W^{\epsilon}((\mathcal{M}_L)\vert_M)$. It follows that $${\rm \overline{Herm}}^{\epsilon}(\mathcal{F}')({\rm
\overline{Herm}}^{\epsilon}(\mathcal{R}^M_{L/K})(N, g))={\rm \overline{Herm}}^{\epsilon}(\mathcal{F}')({\rm
\overline{Herm}}^{\epsilon}(\mathcal{R}^M_{L/K})(N', g'))$$ in $W^{\epsilon \epsilon_0}(E \otimes_K L, \sigma \otimes {\rm
id}_L)$. As the last diagram above commutes, we obtain 
$${\rm \overline{Herm}}^{\epsilon \epsilon_0}(R_{L/K})({\rm \overline{Herm}}^{\epsilon}(\mathcal{F})(N, g))={\rm \overline{Herm}}^{\epsilon \epsilon_0}(R_{L/K})({\rm
\overline{Herm}}^{\epsilon}(\mathcal{F})(N', g')).$$
By [4], proposition 1.2, the map ${\rm \overline{Herm}}^{\epsilon \epsilon_0}(R_{L/K})$ is injective, so we deduce that ${\rm
\overline{Herm}}^{\epsilon}(\mathcal{F})(N, g) = {\rm \overline{Herm}}^{\epsilon}(\mathcal{F})(N', g')$ in $W^{\epsilon \epsilon_0}(E,
\sigma)$.
Since in the category $\mathcal{M}\vert_M$ all idempotents split, the functor $\mathcal{F}$ is an
equivalence of categories and hence the group homomorphism 
${\rm \overline{Herm}}^{\epsilon}(\mathcal{F}):
W^{\epsilon}(\mathcal{M}\vert_M) \rightarrow W^{\epsilon \epsilon_0}(E, \sigma)$ is bijective (see \S 2.1). It follows that
$(N, g)=(N', g')$ in $W^{\epsilon}(\mathcal{M}\vert_M)$. We conclude that the map ${\rm
\overline{Herm}}^{\epsilon}(\mathcal{R}^M_{L/K})$ is injective. 

\medskip

\textbf{Proof of theorem 4.1 } Let $(N,g) \in \mathcal{H}^{\epsilon}(\mathcal{M})$ be such that ${\rm \overline{Herm}}^{\epsilon}(\mathcal{R}_{L/K})(N,g)$ is hyperbolic. This means that $(N,g \otimes 1)$ is a 
hyperbolic unimodular $\epsilon$-hermitian form in the category $\mathcal{M}_L$, so in the category $\mathcal{M}_L\vert_N$ too. Hence 
${\rm \overline{Herm}}^{\epsilon}(\mathcal{R}^N_{L/K})(N,g)=0$ and from to the
injectivity of the map ${\rm \overline{Herm}}^{\epsilon}(\mathcal{R}^N_{L/K})$ we deduce that $(N,g)$ is hyperbolic in $\mathcal{M}\vert_N$, so in $\mathcal{M}$ too. In conclusion the map
${\rm \overline{Herm}}^{\epsilon}(\mathcal{R}_{L/K})$ is injective.

\bigskip
\bigskip

{\bf \S 5. Systems of sesquilinear forms}

\medskip

In this section we generalize results proven in [5] to systems of sesquilinear forms. We first construct an equivalence of categories as in [5], \S 4.

\medskip

{\bf \S 5.1. Preliminaries}

\medskip
Let $A$ be a ring with an involution $\sigma$ and $I$ be a set. 
A \textit{system of sesquilinear forms} over $(A, \sigma)$ is $(V, (s_i)_{i \in I})$, where $V$ is a reflexive right $A$-module of finite type and for all $i \in I$, 
$(V, s_i)$ is a sesquilinear form over $(A, \sigma)$. 
A morphism between two systems of sesquilinear forms $(V, (s_i)_{i \in I})$ and $(V', (s'_i)_{i \in I})$ consists of an isomorphism of $A$-modules
$f:V \tilde{\rightarrow} V'$ such that for every $i \in I$ and $x, y \in V$ we have $s'_i(f(x), f(y))=s_i(x,y)$. Let us denote by $\mathcal{S}^{(I)}(A,\sigma)$ the category of systems of sesquilinear
forms over $(A, \sigma)$.

\medskip

Denote by $J$ the disjoint union of two copies of $I$. We define the \textit{category of $J$-arrows between reflexive $A$-modules} as being the category $\mathcal{M}^{(J)}$ constructed
in the following way:
its objects are of the form $(V, W, (f_i, g_i)_{i \in I})$, where $V$ and $W$ are two reflexive $A$-modules of finite type and for all $i \in I$, $f_i, g_i:V \rightarrow W$ are homomorphisms of
$A$-modules. A morphism from $(V, W, (f_i, g_i)_{i \in I})$ to $(V', W', (f'_i, g'_i)_{i \in I})$ is a pair $(\phi, \psi)$, where $\phi:V \tilde{\rightarrow} V'$ and 
$\psi: W \tilde{\rightarrow} W'$ are isomorphisms of $A$-modules
such that for all $i \in I$ we have $\psi f_i=f'_i \phi$ and $\psi g_i=g'_i \phi$. By defining direct sums in the obvious way we see that $\mathcal{M}^{(J)}$ is an additive category. Let
$(W^*, V^*, (g_i^*, f_i^*)_{i \in I})$ be the dual of $(V, W, (f_i, g_i)_{i \in I})$ and set $E_{(V, W, (f_i, g_i)_{i \in I})}=(e_V, e_W)$. This defines a duality on the category $\mathcal{M}^{(J)}$. 
We observe that if the set $I$ has one element, then $\mathcal{M}^{(J)}$ coincides with the category constructed in [5], \S 3.

\medskip 

We define a functor $\Psi: \mathcal{S}^{(I)}(A, \sigma) \rightarrow \mathcal{M}^{(J)}$ in the following way: Let $(V, (s_i)_{i \in I})$ be an object of $\mathcal{S}^{(I)}(A, \sigma)$ and
for all $i \in I$ let $s_{il}: V \rightarrow V^*$ and $s_{ir}: V \rightarrow V^*$ be the left, respectively the right adjoint of $(V, s_i)$ (cf. \S 1). Then
$(V, V^*, (s_{il}, s_{ir})_{i \in I})$ is an object of $\mathcal{M}^{(J)}$ and it is easy to check that $(e_V, {\rm id}_{V^*})$ defines a unimodular hermitian form on it. Let 
$\Psi(V, (s_i)_{i \in I})=((V, V^*, (s_{il}, s_{ir})_{i \in I}), (e_V, {\rm id}_{V^*}))$. For a morphism $\varphi: (V, (s_i)_{i \in I}) \rightarrow (V', (s'_i)_{i \in I})$ of systems of
sesquilinear forms over $(A, \sigma)$ set $\Psi(\varphi)=(\varphi, \varphi^{*-1})$.

\medskip 

\textbf{5.1.1. Theorem } \textit{The functor $\Psi$ is an equivalence of categories between $\mathcal{S}^{(I)}(A, \sigma)$ and $\mathcal{H}(\mathcal{M}^{(J)})$. }

\medskip 

\textbf{Proof. } The proof is analogous to the one of [5], theorem 4.1. The functor $\Phi: \mathcal{H}(\mathcal{M}^{(J)}) \rightarrow \mathcal{S}(A,\sigma)$ which realizes, together with $\Psi$, the
desired equivalence of categories is defined in the following way: Let $((V, W, (f_i, g_i)_{i \in I}), (\varphi, \psi))$ be an object of $\mathcal{H}(\mathcal{M}^{(J)})$. For every $i \in I$ define 
a sesquilinear form $s_i: V \times V \rightarrow A$ by $s_i(x,y)=(\psi g_i)(x)(y)$ for all $x,y \in V$ and set $\Phi((V, W, (f_i, g_i)_{i \in I}), (\varphi, \psi))=(V, (s_i)_{i \in I})$. For a
morphism $\lambda=(\lambda_1, \lambda_2)$ between two objects of $\mathcal{H}(\mathcal{M}^{(J)})$ set $\Phi(\lambda)=\lambda_1$.

\medskip

If $(V, (s_i)_{i \in I})$ is a system of sesquilinear forms over $(A, \sigma)$, then we denote by $q(V, (s_i)_{i \in I})$ the corresponding object
$(V, V^*, (s_{il}, s_{ir})_{i \in I})$ of the category $\mathcal{M}^{(J)}$. We next describe, following [5], \S 5, the set of isometry classes of systems of sesquilinear forms over $(A, \sigma)$ corresponding
by theorem 5.1.1 to a given object of the category $\mathcal{M}^{(J)}$. 

\medskip

Let us fix a system of sesquilinear forms $(V_0, (s_{0i})_{i \in I})$ over $(A, \sigma)$ and consider the unimodular hermitian form $\eta_{V_0}=(e_{V_0}, {\rm id}_{V_0^*})$ on
$Q_0=q(V_0, (s_{0i})_{i \in I})$. Let $E$ be the endomorphism ring of the object $Q_0$ in $\mathcal{M}^{(J)}$. The form $\eta_{V_0}$ induces an involution $\widetilde{~}$ on $E$, defined by
$\widetilde{f}=\eta_{V_0}^{-1} f^* \eta_{V_0}$ for all $f \in E$, where $f^*$ denotes the dual of $f$ in $\mathcal{M}^{(J)}$. As in [5], \S 5 we denote by $H(\widetilde{~}, E^{\times})$
the set of equivalence classes for the equivalence relation defined on $\{f \in E \vert~\widetilde{f}=f\}$ by: $f \equiv f'$ if there exists a $g \in E^{\times}$ such that $\widetilde{g} fg=f'$. 
Analogously
to [5], theorem 5.1 we obtain:

\medskip

\textbf{5.1.2. Theorem } \textit{The set of isometry classes of systems of sesquilinear forms $(V, (s_i)_{i \in I})$ over $(A, \sigma)$ such that $q(V, (s_i)_{i \in I}) \simeq Q_0$ is in
bijection with $H(\widetilde{~}, E^{\times})$.}

\medskip

\textbf{Proof. } The proof is analogous to the one of [5], theorem 5.1. For every system of sesquilinear forms $(V, (s_i)_{i \in I})$ over $(A, \sigma)$ such that $q(V, (s_i)_{i \in I}) \simeq Q_0$
the unimodular hermitian form $(e_V, {\rm id}_{V^*})$ on $q(V, (s_i)_{i \in I})$ induces a unimodular hermitian form $\eta_V$ on $Q_0$. The desired bijection is given by:
$$\{[(V, (s_i)_{i \in I})]~\vert~(V, (s_i)_{i \in I})\in \mathcal{S}^{(I)}(A, \sigma), ~q(V, (s_i)_{i \in I}) \simeq Q_0\} \rightarrow H(\widetilde{~}, E^{\times})$$
$$[(V, (s_i)_{i \in I})] \mapsto \eta_{V_0}^{-1} \eta_{V}.$$

\bigskip
\bigskip

{\bf \S 5.2. Witt's cancellation theorem}
\medskip \\
Let $K$ be a field of characteristic $\neq 2$, $A$ be a finite-dimensional $K$-algebra and $\sigma$ be an involution on $A$. Analogously to [5], theorem 6.1, a cancellation theorem
holds for systems of sesquilinear forms over $(A, \sigma)$:

\medskip

\textbf{5.2.1. Theorem } \textit{Let $(V, (s_i)_{i \in I})$, $(V', (s'_i)_{i \in I})$ and $(V'', (s''_i)_{i \in I})$ be systems of sesquilinear forms over $(A, \sigma)$ such that
$$(V', (s'_i)_{i \in I}) \oplus (V, (s_i)_{i \in I}) \simeq (V'', (s''_i)_{i \in I}) \oplus (V, (s_i)_{i \in I}).$$
Then we have $(V', (s'_i)_{i \in I}) \simeq (V'', (s''_i)_{i \in I})$. 
}

\medskip

Due to the equivalence between the categories $\mathcal{S}^{(I)}(A, \sigma)$ and $\mathcal{H}(\mathcal{M}^{(J)})$ given by theorem 5.1.1, it is enough to prove that Witt's cancellation theorem holds 
in the category $\mathcal{H}(\mathcal{M}^{(J)})$. This can be proven as in [5], proposition 6.2.

\bigskip
\bigskip

{\bf \S 5.3. Springer's theorem}
\medskip \\
In this subsection we prove an analogue of Springer's theorem for systems of sesquilinear forms defined over finite-dimensional algebras with involution.
\medskip \\
Let $K$ be a field of characteristic $\neq 2$, $A$ be a finite-dimensional $K$-algebra and $\sigma$ be a $K$-linear involution on $A$. We also consider a finite extension $L$ of $K$,
the finite-dimensional $L$-algebra $A_L=A \otimes_K L$ and the $L$-linear involution $\sigma_L=\sigma \otimes {\rm id}_L$ on $A_L$. 
If $(V, (s_i)_{i \in I})$ is a system of sesquilinear forms over $(A, \sigma)$, then
we denote by $(V, (s_i)_{i \in I})_L=(V_L, ((s_i)_L)_{i \in I})$ the system of sesquilinear forms over $(A_L, \sigma_L)$ obtained by extension of scalars.
\medskip \\
\textbf{5.3.1. Theorem } \textit{Suppose that $L/K$ is an extension of odd degree. Let $(V, (s_i)_{i \in I})$ and $(V', (s'_i)_{i \in I})$ be two
systems of sesquilinear forms over $(A, \sigma)$. If $(V, (s_i)_{i \in I})_L$ and $(V', (s'_i)_{i \in I})_L$ are isometric over $(A_L, \sigma_L)$, then 
$(V, (s_i)_{i \in I})$ and $(V', (s'_i)_{i \in I})$ are isometric over $(A, \sigma)$. }
\medskip \\
\textbf{Proof. } The proof is analogous to the one of [5], theorem 7.1 and uses the equivalence of categories given by theorem 5.1.1. 

\bigskip
\bigskip

{\bf \S 5.4. Weak Hasse principle}

\medskip
The well-known weak Hasse principle states that if two quadratic forms defined over a global field $k$ of characteristic $\neq 2$ 
become isometric over all the completions of $k$, then they are already isometric over $k$. The aim of this section is to generalize this result to the case
of systems of sesquilinear forms defined over a skew field with involution. Throughout this section we suppose that $I$ is finite.

\medskip

For results concerning the weak Hasse principle for systems of hermitian or quadratic forms 
over fields see [1], respectively [2]. A weak Hasse principle for sesquilinear forms defined over a skew field with involution has been proven in [5]. 

\medskip

Let $K$ be a field of characteristic $\neq 2$, $D$ be a finite-dimensional skew field with center $K$ and $\sigma$ be an involution on $D$. Denote by $k$
the fixed field of $\sigma$ in $K$. Then either $k=K$ (when $\sigma$ is said to be of the \textit{first kind}) or $K$ is a quadratic extension of $k$ and the restriction of 
$\sigma$ to $K$ is the non-trivial automorphism of $K$ over $k$ (in which case $\sigma$ is said to be of the \textit{second kind} or a \textit{unitary involution}). 

\medskip

Suppose that $k$ is a global field. For every prime spot $p$ of $k$, let $k_p$ be the completion of $k$ at $p$, $K_p=K \otimes_k k_p$ and $D_p=D \otimes_k k_p$. Then $D_p$ is an algebra with center
$K_p$ and consider on it the involution $\sigma_p=\sigma \otimes {\rm id}_{k_p}$. Then from any sesquilinear form $(V,s)$ over $(D, \sigma)$ we obtain by extension of scalars a sesquilinear form
$(V_p, s_p)$ over $(D_p, \sigma_p)$, where $V_p=V \otimes_k k_p$ is a free right $D_p$--module of rank ${\rm dim}_D(V)$. Hence we obtain a notion of extension of scalars for systems of
sesquilinear forms over $(D, \sigma)$ by setting $(V, (s_i)_{i \in I})_p=(V_p, ((s_i)_p)_{i \in I})$.

\medskip

We say that the \textit{weak Hasse principle} holds for systems of sesquilinear
forms over $(D, \sigma)$ if any two systems of sesquilinear forms $(V,(s_i)_{i \in I})$ and $(V',(s'_i)_{i \in I})$ over $(D, \sigma)$ that become isometric over all the completions of $k$ 
(i.e. $(V_p, ((s_i)_p)_{i \in I}) \simeq (V'_p, ((s'_i)_p)_{i \in I})$ over
$(D_p, \sigma_p)$ for every prime spot $p$ of $k$) are already isometric over $(D, \sigma)$. 

\medskip

For the problem of determining when the weak Hasse principle holds we will use the ideas developed in [5].
We denote by $\mathcal{M}^{(J)}$ ( $(\mathcal{M}_p)^{(J)}$ ) the category of $J$-arrows between reflexive $D$-- (resp. $D_p$--) modules. Let us fix a system of sesquilinear forms 
$s_0$ over $(D, \sigma)$.
Consider a system of sesquilinear forms $s$ over $(D, \sigma)$ that becomes isometric to $s_0$ over all the completions of $k$. Denote by $q(s)$ the object of $\mathcal{M}^{(J)}$ corresponding
to $s$ by the equivalence of categories given in theorem 5.1.1. For every prime spot $p$ of $k$ we have
$(s_0)_p \simeq s_p$, hence $q((s_0)_p) \simeq q(s_p)$. Since $q$ commutes with base change, we obtain $q(s_0)_p \simeq q(s)_p$
and since $I$ is finite, it follows that $q(s_0) \simeq q(s)$. 

\medskip 

Denote by $E$ the endomorphism ring of $q(s_0)$ in $\mathcal{M}$. The unimodular hermitian form $\eta_{s_0}=(e_{V_0}, {\rm id}_{V_0^*})$ defined on $q(s_0)$ induces an 
involution $\widetilde{~}$ on $E$. For every $p$ let $E_p=E \otimes_k k_p$, on which we consider the involution $\widetilde{~}_p=\widetilde{~} \otimes {\rm id}_{k_p}$. By theorem 5.1.2 it is clear that if the localisation map
$$\Phi: H(\widetilde{~}, E^{\times}) \rightarrow \prod_p H(\widetilde{~}_p, E_p^{\times})$$
is injective, then the weak Hasse principle holds for systems of sesquilinear forms over $(D, \sigma)$. As in [5], \S 8, we prove the following:

\medskip

\textbf{5.4.1. Theorem } \textit{If $\sigma$ is a unitary involution, then the map $\Phi$ is injective and hence the weak Hasse principle holds for systems of sesquilinear forms over $(D, \sigma)$. }

\bigskip
\bigskip

{\bf \S 5.5. Finiteness results }

\medskip

In this section we generalize theorem 3.3.1 and the results proven in [5], \S 9 to systems of sesquilinear forms. Following \S 2.3 we say that two systems of
sesquilinear forms over $(A, \sigma)$ are in the same genus if they become isometric over all the extensions of $A$ to $\mathbb{Z}_p$, where $p$ is a prime. 

\medskip

Fix a system of sesquilinear forms $(V, (s_i)_{i \in I})$ over $(A, \sigma)$ and denote by $q(V,(s_i)_{i \in I})$ the corresponding object $(V, V^*, (s_{il}, s_{ir})_{i \in I})$ of $\mathcal{M}^{(J)}$ 
and by $E$ its endomorphism ring in $\mathcal{M}^{(J)}$. It is straightforward to prove the following results using the methods of [5], \S 9 and the equivalence of categories 
$F: \mathcal{S}^{(I)}(A, \sigma) \rightarrow \mathcal{H}(\mathcal{M}^{(J)})$ given by theorem 5.1.1: 

\medskip 

\textbf{5.5.1. Theorem } \textit{If there exists a non-zero integer $m$ such that ${\rm End}_A(V)$ is $\mathbb{Z}[1/m]$-finite, then there exist only finitely many isometry classes of 
systems of sesquilinear forms on $V$.}

\medskip

\textbf{5.5.2. Theorem } \textit{Let $N$ be an object of $\mathcal{M}^{(J)}$ and assume that there exists a non-zero integer $m$ such that ${\rm End}_{\mathcal{M}^{(J)}}(N)$ is $\mathbb{Z}[1/m]$-finite.
Then there exist only finitely many isometry classes of systems of sesquilinear forms $(V, (s_i)_{i \in I})$ over $(A, \sigma)$ such that $q(V, (s_i)_{i \in I}) \simeq N$.}

\medskip

\textbf{5.5.3. Theorem } \textit{If there exists a non-zero integer $m$ such that ${\rm End}_A(V)$ is $\mathbb{Z}[1/m]$-finite, then $(V, (s_i)_{i \in I})$ contains only finitely many isometry
classes of orthogonal summands.}

\medskip

\textbf{5.5.4. Theorem } \textit{If ${\rm End}_A(V)$ is $\mathbb{Q}$-finite, then the genus of $(V, (s_i)_{i \in I})$ contains only a finite number of isometry classes of
systems of sesquilinear forms.}

\bigskip
\bigskip

\textbf{Aknowledgements } The second named author would like to thank Emmanuel Lequeu for many interesting and useful discussions.

\newpage

{\bf Bibliography}
\medskip \\

\noindent
[1] E. Bayer--Fluckiger, Principe de Hasse faible pour les syst\`emes de formes quadratiques, J. reine angew. Math. 378 (1987), 53-59. \\

\noindent
[2] E. Bayer--Fluckiger, Intersection de groupes orthogonaux et principe de Hasse faible pour les syst\`emes de formes quadratiques sur un corps global, 
C.R. Acad. Sci. Paris I, 301 (1985), 911-914. \\

\noindent
[3] E. Bayer--Fluckiger, C. Kearton, S.M. Wilson, Hermitian forms in additive categories: finiteness results, 
{\it J. Algebra} {\bf 123} (1989), no. 2, 336--350. \\

\noindent
[4] E. Bayer--Fluckiger, H.W. Lenstra, Jr.,  Forms in odd degree
extensions and self-dual normal bases,
{\it Amer. J. Math.} {\bf 112} (1990), 359-373. \\

\noindent
[5] E. Bayer--Fluckiger, D.A. Moldovan, Sesquilinear forms over rings with involution, preprint \\



\noindent
[6] M. Knus, {\it Quadratic and hermitian forms over rings}, Grundlehren der Math. Wiss.
{\bf 294}, Springer--Verlag (1991). \\

\noindent
[7] H.--G. Quebbemann, W. Scharlau, M. Schulte, Quadratic and hermitian
forms in additive and abelian categories, {\it J. Algebra} {\bf 59} (1979), 264--289. \\

\noindent
[8] W. Scharlau, {\it Quadratic and Hermitian Forms}, Grundlehren der Math. Wiss. {\bf 270},
Springer--Verlag (1985).

\vfill
\eject
\end{document}